\documentclass[12pt]{article}

\usepackage{amsmath, epsfig, cite}
\usepackage{amssymb}
\usepackage{amsfonts}
\usepackage{latexsym}
\usepackage{float}
\usepackage{color}

\newtheorem{thm}{Theorem}[section]

\newtheorem{cor}[thm]{Corollary}
\newtheorem{conj}[thm]{Conjecture}
\newtheorem{lem}[thm]{Lemma}

\numberwithin{equation}{section}

\newcommand{\qed}{{\hfill$\square$}\medskip}

\begin{document}

%\linenumbers
\begin{center}
{\large\bf On two supercongruences of double binomial sums}
\end{center}

\vskip 2mm \centerline{Long Li$^1$ and Ji-Cai Liu$^2$\footnote{Corresponding author.} }
\begin{center}
{\footnotesize  $^1$School of Mathematics and Statistics, Huaiyin Normal
University, Huai'an 223300, PR China\\
{\tt lli@hytc.edu.cn } \\[10pt]
$^2$Department of Mathematics, Wenzhou University, Wenzhou 325035, PR China\\
{\tt jcliu2016@gmail.com }}
\end{center}

%%date: November 27, 2014
%%date : December 4, 2014

\vskip 0.7cm \noindent{\bf Abstract.}
In this note, we confirm two conjectural supercongruences on double sums of
binomial coefficients due to El Bachraoui.


\vskip 2mm
\noindent{\it MR Subject Classifications}: 11A07, 11B65, 05A19

\section{Introduction}
Recall that the $q$-shifted factorials are given by $(a;q)_n=(1-a)(1-aq)\cdots(1-aq^{n-1})$ for $n\ge 1$ and $(a;q)_0=1$, and the $q$-integers are defined by $[n]=(1-q^n)/(1-q)$.
For polynomials $A_1(q), A_2(q),P(q)\in \mathbb{Z}[q]$, the $q$-congruence $$A_1(q)/A_2(q)\equiv 0\pmod{P(q)}$$ is understood as $A_1(q)$ is divisible by $P(q)$ and $A_2(q)$ is coprime with $P(q)$. In general, for rational functions $A(q),B(q)\in \mathbb{Z}(q)$,
\begin{align*}
A(q)\equiv B(q)\pmod{P(q)}\Longleftrightarrow
A(q)-B(q)\equiv 0\pmod{P(q)}.
\end{align*}

In the past few years, $q$-congruence for sums of binomial coefficients as well as hypergeometric series attracted many experts' attention (see, for instance, \cite{bachraoui-rj-2019,guo-rm-2020,gl-jdea-2018,gs-rm-2019,gz-am-2019,lh-bams-2020,lp-aam-2020}). In particular, Guo and Zudilin \cite{gz-am-2019} developed a creative microscoping method to prove many interesting $q$-congruences, such as
\begin{align}
\sum_{k=0}^{n-1}(-1)^k\frac{(q;q^2)_k(-q;q^2)_k^2}{(q^4;q^4)_k(-q^4;q^4)_k^2}
[6k+1]q^{3k^2}\equiv 0\pmod{[n]},\label{new-1}
\end{align}
and
\begin{align}
\sum_{k=0}^{n-1}\frac{(q^2;q^4)_k(-q;q^2)_k^2}{(q^4;q^4)_k(-q^4;q^4)_k^2}
[6k+1]q^{k^2}\equiv 0\pmod{[n]},\label{new-2}
\end{align}
for any odd positive integer $n$.

Motivated by \eqref{new-1} and \eqref{new-2}, El Bachraoui \cite[Theorems 1 and 2]{bachraoui-rj-2019} established the following two $q$-congruences:
\begin{align}
\sum_{k=0}^{n-1}\sum_{j=0}^k c_q(j)c_q(k-j)\equiv 0\pmod{[n]},\label{new-3}
\end{align}
and
\begin{align}
\sum_{k=0}^{n-1}\sum_{j=0}^k c'_q(j)c'_q(k-j)\equiv 0\pmod{[n]},\label{new-4}
\end{align}
where $c_q(k)$ and $c'_q(k)$ denote the $k$-th term of the summations on the left-hand sides of
\eqref{new-1} and \eqref{new-2}, respectively.

Suppose $p$ is an odd prime.
Letting $q\to 1$ and $n=p$ in \eqref{new-3} and \eqref{new-4} gives
\begin{align}
\sum_{k=0}^{p-1}\left(-\frac{1}{8}\right)^k\sum_{j=0}^k{2j\choose j}{2k-2j\choose k-j}(6j+1)(6k-6j+1)\equiv 0\pmod{p},\label{new-5}
\end{align}
and
\begin{align}
\sum_{k=0}^{p-1}\left(\frac{1}{4}\right)^k\sum_{j=0}^k{2j\choose j}{2k-2j\choose k-j}(6j+1)(6k-6j+1)\equiv 0\pmod{p}.\label{new-6}
\end{align}

El Bachraoui \cite[Conjectures 1 and 2]{bachraoui-rj-2019} also conjectured two extensions of
\eqref{new-5} and \eqref{new-6} as follows:
\begin{conj}
For any odd prime $p$, we have
\begin{align}
\sum_{k=0}^{p-1}\left(-\frac{1}{8}\right)^k\sum_{j=0}^k{2j\choose j}{2k-2j\choose k-j}(6j+1)(6k-6j+1)\equiv -\frac{p}{2}\pmod{p^2},\label{a-1}
\end{align}
and
\begin{align}
\sum_{k=0}^{p-1}\left(\frac{1}{4}\right)^k\sum_{j=0}^k{2j\choose j}{2k-2j\choose k-j}(6j+1)(6k-6j+1)\equiv p\pmod{p^2}.
\label{a-2}
\end{align}
\end{conj}

The above two conjectural supercongruences motivate us to establish the following more general result, which includes \eqref{a-1} and \eqref{a-2} as special cases.

\begin{thm}\label{t-1}
For any positive integer $n$, we have
\begin{align}
&\sum_{k=0}^{n-1}\left(\frac{q}{4}\right)^k\sum_{j=0}^k{2j\choose j}{2k-2j\choose k-j}(6j+1)(6k-6j+1)\notag\\[7pt]
&=\frac{(9n^2-15n+8)q^{n+2}-(18n^2-12n-8)q^{n+1}+(9n^2+3n+2)q^n-2(2q+1)^2}{2(1-q)^3}.\label{a-3}
\end{align}
\end{thm}

Letting $q\to -\frac{1}{2}$ and $q\to 1$ in \eqref{a-3}, we obtain the following two combinatorial identities:
\begin{cor}
For any positive integer $n$, we have
\begin{align}
\sum_{k=0}^{n-1}\left(-\frac{1}{8}\right)^k\sum_{j=0}^k{2j\choose j}{2k-2j\choose k-j}(6j+1)(6k-6j+1)=\left(-\frac{1}{2}\right)^n n(1-3n),\label{a-4}
\end{align}
and
\begin{align}
\sum_{k=0}^{n-1}\left(\frac{1}{4}\right)^k\sum_{j=0}^k{2j\choose j}{2k-2j\choose k-j}(6j+1)(6k-6j+1)=\frac{n(3n^2-3n+2)}{2}.\label{a-5}
\end{align}
\end{cor}

It is clear that \eqref{a-1} and \eqref{a-2} can be deduced from \eqref{a-4} and \eqref{a-5} directly. We shall prove Theorem \ref{t-1} in the next section.

\section{Proof of Theorem \ref{t-1}}
In order to prove Theorem \ref{t-1}, we need the following preliminary result.
\begin{lem}
For any non-negative integer $n$, we have
\begin{align}
\sum_{j=0}^k{2j\choose j}{2k-2j\choose k-j}(6j+1)(6k-6j+1)=
4^k\left(\frac{9}{2}k^2+\frac{3}{2}k+1\right).\label{b-1}
\end{align}
\end{lem}
{\it Proof.}
Note that central binomial coefficients possess the following generating function:
\begin{align}
\frac{1}{\sqrt{1-4x}}=\sum_{j=0}^{\infty}{2j\choose j}x^j.\label{d-1}
\end{align}
On the other hand,
\begin{align*}
\frac{1}{1-4x}=\sum_{k=0}^{\infty}(4x)^k.
\end{align*}
Thus,
\begin{align}
\left(\sum_{j=0}^{\infty}{2j\choose j}x^j\right)^2=\sum_{k=0}^{\infty}(4x)^k.\label{c-1}
\end{align}
Comparing the coefficient of $x^k$ on both sides of \eqref{c-1}, we obtain
\begin{align}
\sum_{j=0}^k{2j\choose j}{2k-2j\choose k-j}=4^k.\label{c-2}
\end{align}

Differentiating both sides of \eqref{d-1} respect to $x$, we obtain
\begin{align}
\frac{2x}{\sqrt{(1-4x)^3}}=\sum_{j=0}^{\infty}j{2j\choose j}x^j.\label{d-2}
\end{align}
On the other hand,
\begin{align*}
\frac{4x^2}{(1-4x)^3}=4x^2\sum_{k=0}^{\infty}{-3\choose k}(-4x)^k.
\end{align*}
Since
\begin{align*}
{-3\choose k}=\frac{(-1)^k(k+1)(k+2)}{2},
\end{align*}
we have
\begin{align}
\frac{4x^2}{(1-4x)^3}=\sum_{k=0}^{\infty}\frac{4^{k+1}(k+1)(k+2)}{2}x^{k+2}.\label{d-3}
\end{align}
Noting \eqref{d-2} and \eqref{d-3}, we find that
\begin{align}
\left(\sum_{j=0}^{\infty}j{2j\choose j}x^j\right)^2=\sum_{k=0}^{\infty}\frac{4^{k+1}(k+1)(k+2)}{2}x^{k+2}.\label{d-4}
\end{align}
Comparing the coefficient of $x^k$ on both sides of \eqref{d-4}, we obtain
\begin{align}
\sum_{j=0}^k{2j\choose j}{2k-2j\choose k-j}j(k-j)=\frac{4^kk(k-1)}{8}.\label{d-5}
\end{align}

Finally, using \eqref{c-2} and \eqref{d-5} we arrive at
\begin{align*}
&\sum_{j=0}^k{2j\choose j}{2k-2j\choose k-j}(6j+1)(6k-6j+1)\\
&=(6k+1)\sum_{j=0}^k{2j\choose j}{2k-2j\choose k-j}+36\sum_{j=0}^k{2j\choose j}{2k-2j\choose k-j}j(k-j)\\
&=4^k\left(\frac{9}{2}k^2+\frac{3}{2}k+1\right),
\end{align*}
as desired.
\qed

{\noindent \it Proof of Theorem \ref{t-1}.}
By \eqref{b-1}, we have
\begin{align}
&\sum_{k=0}^{n-1}\left(\frac{q}{4}\right)^k\sum_{j=0}^k{2j\choose j}{2k-2j\choose k-j}(6j+1)(6k-6j+1)\notag\\
&=\sum_{k=0}^{n-1}q^k\left(\frac{9}{2}k^2+\frac{3}{2}k+1\right).\label{b-2}
\end{align}
Let
\begin{align*}
&S_n=\sum_{k=0}^{n-1} q^k k\quad \text{and}\quad T_n=\sum_{k=0}^{n-1} q^k k^2.
\end{align*}
Note that
\begin{align*}
(1-q)S_{n}
&=\sum_{k=0}^{n-1} q^k k-\sum_{k=0}^{n-1} q^{k+1} k\\
&=\sum_{k=1}^{n-1} q^k k-\sum_{k=1}^{n}q^{k} (k-1)\\
&=\sum_{k=1}^{n-1} q^k k-\sum_{k=1}^{n-1} q^{k} (k-1)-q^{n}(n-1)\\
&=\sum_{k=1}^{n-1} q^k-q^{n}(n-1).
\end{align*}
Thus,
\begin{align}
S_n=\frac{q(1-q^{n-1})}{(1-q)^2}-\frac{q^n(n-1)}{1-q}.\label{b-3}
\end{align}

On the other hand,
\begin{align}
(1-q)T_{n}
&=\sum_{k=0}^{n-1} q^k k^2-\sum_{k=0}^{n-1} q^{k+1} k^2\notag\\
&=\sum_{k=1}^{n-1} q^k k^2-\sum_{k=1}^{n} q^{k} (k-1)^2\notag\\
&=\sum_{k=1}^{n-1} q^k k^2-\sum_{k=1}^{n-1} q^{k} (k-1)^2-q^{n}(n-1)^2\notag\\
&=\sum_{k=1}^{n-1} q^k(2k-1)-q^{n}(n-1)^2.\label{b-4}
\end{align}
Combining \eqref{b-3} and \eqref{b-4}, we obtain
\begin{align}
T_n=\frac{2q(1-q^{n-1})}{(1-q)^3}-\frac{2q^n(n-1)}{(1-q)^2}-\frac{q(1-q^{n-1})}{(1-q)^2}-\frac{q^{n}(n-1)^2}{1-q}.
\label{b-5}
\end{align}

Finally, substituting \eqref{b-3} and \eqref{b-5} into \eqref{b-2}, we arrive at
\begin{align*}
&\sum_{k=0}^{n-1}\left(\frac{q}{4}\right)^k\sum_{j=0}^k{2j\choose j}{2k-2j\choose k-j}(6j+1)(6k-6j+1)\notag\\[7pt]
&=\frac{(9n^2-15n+8)q^{n+2}-(18n^2-12n-8)q^{n+1}+(9n^2+3n+2)q^n-2(2q+1)^2}{2(1-q)^3},
\end{align*}
as desired.
\qed

\vskip 5mm \noindent{\bf Acknowledgments.}
The first author was supported by the Natural Science Foundation of the Jiangsu Higher Education Institutions of China (grant 19KJB110006).
The second author was supported by the National Natural Science Foundation of China (grant 11801417).


\begin{thebibliography}{99}

\small \setlength{\itemsep}{-.8mm}

\bibitem{bachraoui-rj-2019}M. El Bachraoui, On supercongruences for truncated sums of squares of basic hypergeometric series, Ramanujan J., online, doi: 10.1007/s11139-019-00226-0.

\bibitem{guo-rm-2020}V.J.W. Guo, Proof of some $q$-supercongruences modulo the fourth power of a cyclotomic polynomial, Results Math. 75 (2020), Art. 77.

\bibitem{gl-jdea-2018}V.J.W. Guo and J.-C. Liu, $q$-Analogues of two Ramanujan-type formulas for $1/\pi$, J. Difference Equ. Appl. 24 (2018), 1368--1373.

\bibitem{gs-rm-2019}V.J.W. Guo and M. J. Schlosser, Some new $q$-congruences for truncated basic hypergeometric series: even powers, Results Math. 75 (2020), Art. 1.

\bibitem{gz-am-2019}V.J.W. Guo and W. Zudilin, A $q$-microscope for supercongruences, Adv. Math. 346 (2019), 329--358.

\bibitem{lh-bams-2020}J.-C. Liu and Z.-Y. Huang, A truncated identity of Euler and related $q$-congruences, Bull. Aust. Math. Soc., online, doi:10.1017/S0004972720000301.

\bibitem{lp-aam-2020}J.-C. Liu and F. Petrov, Congruences on sums of $q$-binomial coefficients, Adv. in Appl. Math. 116 (2020), 102003.

\end{thebibliography}
\end{document}